\documentclass[12pt,twoside]{article}
\usepackage{amsthm}
  \usepackage{amsmath}

  \def\im{{\rm Im}}
  \def\Hom{{\rm Hom}}

\def\Ggldim{{\rm Ggldim}}

\def\Ext{{\rm Ext}}

\def\pd{{\rm pd}}
\def\id{{\rm id}}
\def\fd{{\rm fd}}

\def\Gpd{{\rm Gpd}}
\def\Gid{{\rm Gid}}

\def\gd{{\rm Gdim}}

\date{}

\def\proof{{\noindent\sl Proof.~~}}
\begin{document}

\centerline{}

\centerline{}

\centerline {\Large{\bf On $n$-strongly Gorenstein rings}}

\centerline{}

\centerline{\bf {Mohamed Chhiti,  Khalid Louartiti and Mohammed
Tamekkante }}

\centerline{}

\centerline{Department of Mathematics, }

\centerline{Faculty of Science and Technology of Fez, Box 2202}

\centerline{University S.M. Ben Abdellah Fez, }

\centerline{Morocco}

\centerline{}

\newtheorem{Theorem}{\quad Theorem}[section]

\newtheorem{Definition}[Theorem]{\quad Definition}
\newtheorem{Definitions}[Theorem]{\quad Definitions}
\newtheorem{Corollary}[Theorem]{\quad Corollary}
\newtheorem{Proposition}[Theorem]{\quad Proposition}

\newtheorem{Lemma}[Theorem]{\quad Lemma}

\newtheorem{Example}[Theorem]{\quad Example}
\newtheorem{Remark}[Theorem]{\quad Remark}
\begin{abstract} This paper introduces and studies a particular subclass of the
class  of commutative rings with finite Gorenstein global
dimension.

\end{abstract}

{\bf Mathematics Subject Classification:} 13D05, 13D02 \\

{\bf Keywords:} strongly (n-)Gorenstein projective and injective
modules, Gorenstein global dimensions.

\section{Introduction}
Throughout the paper, all rings are commutative with identity, and
all modules are unitary.\\ Let $R$ be a ring, and let $M$ be an
$R$-module. As usual, we use $\pd_R(M)$, $\id_R(M)$, and
$\fd_R(M)$ to denote, respectively, the classical projective
dimension, injective dimension, and flat
dimension of $M$.\\

For a two-sided Noetherian ring $R$, Auslander and Bridger
\cite{Aus bri} introduced the $G$-dimension, $\gd_R (M)$, for
every finitely generated $R$-module $M$. They showed that  $\gd_R
(M)\leq \pd_R (M)$ for all finitely generated $R$-modules $M$, and
equality holds if $\pd_R (M)$ is finite.

Several decades later, Enochs and Jenda \cite{Enochs,Enochs2}
introduced the notion of Gorenstein projective dimension
($G$-projective dimension for short), as an extension of
$G$-dimension to modules that are not necessarily finitely
generated, and the Gorenstein injective dimension ($G$-injective
dimension for short) as a dual notion of Gorenstein projective
dimension. Then, to complete the analogy with the classical
homological dimension, Enochs, Jenda, and Torrecillas \cite{Eno
Jenda Torrecillas} introduced the Gorenstein flat dimension. Some
references are
 \cite{Bennis and Mahdou2, Christensen, Christensen
and Frankild, Enochs, Enochs2, Eno Jenda Torrecillas, Holm}.\\

Recall that an $R$-module $M$ is called Gorenstein projective, if
there exists an exact sequence of projective  $R$-modules:
$$\mathbf{P}:...\rightarrow P_1\rightarrow P_0\rightarrow
P^0\rightarrow P^1\rightarrow ...$$ such that $M\cong
\im(P_0\rightarrow P^0)$ and such that the functor $\Hom_R(-,Q)$
leaves $\mathbf{P}$ exact whenever $Q$ is a projective $R$-module.
The complex $\mathbf{P}$ is called a complete projective
resolution. The  Gorenstein injective $R$-modules are defined
dually.\\
The  Gorenstein projective and  injective
dimensions are defined in terms of resolutions and denoted by $\Gpd(-)$ and  $\Gid(-)$, respectively (\cite{Christensen, Enocks and janda, Holm}).\\

In \cite{Bennis and Mahdou2}, the authors proved, for any
associative ring $R$, the equality
$$\sup\{\Gpd_R(M)\mid \text{M is a (left)  R-module}\}=\sup\{\Gid_R(M)\mid \text{M is a (left)  R-module}\}.$$
They called the common value of the above quantities the left
Gorenstein global dimension of $R$ and denoted it by
$l.\Ggldim(R)$.  Since in this paper all rings are commutative, we
drop the letter $l$.

Recently, in \cite{MahdouTamekkante}, particular modules of finite
Gorenstein projective, injective, and flat dimensions are defined
as follows:
\begin{Definitions}
Let $n$ be a positive integer.
\begin{enumerate}
  \item An $R$-module $M$ is said to be strongly $n$-Gorenstein projective,
if there exists a short exact sequence of $R$-modules
$0\longrightarrow M \longrightarrow P \longrightarrow M
\longrightarrow 0$ where $\pd_R(P)\leq n$ and
$\Ext^{n+1}_R(M,Q)=0$ whenever $Q$ is projective.
  \item An $R$-module $M$ is said to be strongly $n$-Gorenstein injective,
if there exists a short exact sequence of $R$-modules
$0\longrightarrow M \longrightarrow I \longrightarrow M
\longrightarrow 0$ where $\id_R(I)\leq n$ and
$\Ext^{n+1}_R(E,M)=0$ whenever $E$ is injective.
\end{enumerate}
\end{Definitions}

Clearly, strongly $0$-Gorenstein projective and  injective are
just the strongly Gorenstein projective, injective, and flat
modules, respectively (\cite[Propositions 2.9 and 3.6]{Bennis and
Mahdou2}).\\

In this paper, we investigate these modules to characterize a new
class of rings with finite  Gorenstein global dimension, which we
call $n$-strongly Gorenstein rings.

\section{$n$-Strongly Gorenstein rings}

In \cite{MahdouTamekkante}, the authors proved the following
proposition:

\begin{Proposition}[Proposition 2.16, \cite{MahdouTamekkante}]\label{propo defini}
Let $R$ be a ring. The following statements are equivalent:
\begin{enumerate}
  \item Every module is strongly $n$-Gorenstein projective.
  \item Every module is strongly $n$-Gorenstein injective.
\end{enumerate}
\end{Proposition}

Thus, we give the following definition:
\begin{Definition}Let $n$ be a positive integer. A ring $R$ is called $n$-strongly Gorenstein ($n$-SG ring for
short), if R satisfies one of the equivalent conditions of
Proposition \ref{propo defini}.
\end{Definition}

The  $0$-SG rings and   $1$-SG rings  are already studied in
\cite{ouarghi, MT1} and  they are called strongly Gorenstein
semi-simple rings and  strongly Gorenstein hereditary rings,
respectively. Clearly, by definition, every $n$-SG ring is $m$-SG
  whenever $n\leq m$.

Our first result gives a characterization of strongly
$n$-Gorenstein rings.
\begin{Proposition}\label{prop1} For a ring $R$ and a positive integer $n$, the following statements are equivalent:
\begin{enumerate}
    \item $R$ is an $n$-SG ring.
    \item $\Ggldim(R)\leq n$ and for every $R$-module $M$ there
    exists a short exact sequence of $R$-modules $$0\longrightarrow M \longrightarrow
    P\longrightarrow M \longrightarrow 0$$ where $\pd_R(P)<\infty$.
    \item $\Ggldim(R)<\infty$ and for every $R$-module $M$ there
    exists a short exact sequence of $R$-modules $$0\longrightarrow M \longrightarrow
    P\longrightarrow M \longrightarrow 0$$ where $\pd_R(P)\leq n$.
\end{enumerate}
\end{Proposition}
\proof $(1\Rightarrow 2)$ Clear since for every $n$-SG ring $R$ we
have $\Ggldim(R)\leq n$ (by \cite[Proposition
2.2(1)]{MahdouTamekkante}).

 $(2\Rightarrow 3)$ Follows directly
from \cite[Corollary 2.7]{Bennis and Mahdou2}.

 $(3\Rightarrow 1)$
Follows from \cite[Proposition 2.10]{MahdouTamekkante}.\\

The next result studies the direct product of $n$-SG rings.

\begin{Theorem}\label{product}
Let $\displaystyle\{R_i\}_{i=1}^m$ be a family of rings and set
$R:=\displaystyle\prod_{i=1}^m R_i$. Then, $R$ is an $n$-SG ring
if and only if $R_i$ is an $n$-SG
 ring  for each $i=1,..,m$.
\end{Theorem}
\proof By induction on $m$ it suffices to prove the assertion for
$m=2$. First suppose that $R_1\times R_2$ is an $n$-SG ring. We
claim that $R_1$ is an $n$-SG ring. Let $M$ be an arbitrary $R_1$
module. $M\times 0$ can be viewed as an $R_1\times R_2$-module.
For such module and since $R_1\times R_2$ is an $n$-SG ring, there
is an exact sequence $0\longrightarrow M\times 0\longrightarrow P
\longrightarrow M\times 0\longrightarrow 0$ where $\pd_{R_1\times
R_2}(P)\leq n$. Thus, since $R_1$ is a projective $R_1\times R_2$
module, by applying $-\otimes_{R_1\times R_2}R_1$ to the sequence
above, we find the short exact sequence of $R$-modules:
$0\longrightarrow M\times 0\otimes_{R_1\times
R_2}R_1\longrightarrow P\otimes_{R_1\times R_2}R_1\longrightarrow
M\times 0\otimes_{R_1\times R_2}R_1\longrightarrow 0$. Clearly
$\pd_{R_1}(P\otimes_{R_1\times R_2}R_1)\leq \pd_{R_1\times
R_2}(P)\leq n$. Moreover, we have the isomorphism of $R$-modules:
$$M\times 0\otimes_{R_1\times R_2}R_1\cong M\times
0\otimes_{R_1\times R_2}(R_1\times R_2)/(0\times R_2)\cong M.$$
Thus, we obtain an exact sequence of $R$-module with the form:
$0\longrightarrow M\longrightarrow P\otimes_{R_1\times
R_2}R_1\longrightarrow M\longrightarrow 0$. On the other hand, by
\cite[Theorem 3.1]{Bennis and Mahdou3}, we have $\Ggldim(R_1)\leq
\Ggldim(R_1\times R_2)\leq n$. Thus, using Proposition
\ref{prop1}, $R_1$ is an $n$-SG ring, as desired. By the same
argument, $R_2$ is also an $n$-SG ring.\\ Now, suppose that $R_1$
and $R_2$ are an $n$-SG rings and we claim that $R_1\times R_2$ is
an $n$-SG ring. Let $M$ be an $R_1\times R_2$-module. We have
$$M\cong M\otimes_{R_1\times R_2}(R_1\times R_2)\cong
M\otimes_{R_1\times R_2}((R_1\times 0)\oplus (R_2\times 0))\cong
M_1\times M_2$$ where $M_i=M\otimes_{R_1\times R_2}R_i$ for
$i=1,2$. For each $i=1,2$, there is an exact sequence
$0\longrightarrow M_i\longrightarrow P_i\longrightarrow
M_i\longrightarrow 0$ where $\pd_{R_i}(P_i)\leq n$ since $R_i$ is
an $n$-SG ring. Thus, we have the exact sequence of $R_1\times
R_2$-modules: $$0\longrightarrow M_1\times M_2\longrightarrow
P_1\times P_2\longrightarrow M_1\times M_2\longrightarrow 0.$$ On
the other hand, $\pd_{R_1\times R_2}(P_1\times
P_2)=\sup\{\pd_{R_i}(P_i)\}_{1,2}\leq n$ (by \cite[Lemma 2.5
(2)]{Mahdou costa}). Moreover, by \cite[Theorem 3.1]{Bennis and
Mahdou3}, $\Ggldim(R_1\times R_2)=\sup\{\Ggldim(R_i)\}_{1,2}\leq
n$. Thus, from Proposition \ref{prop1}, $R_1\times R_2$ is an
$n$-SG ring, as desired.\\

Let $T:=R[X_1,X_2,...,X_n]$ be the polynomial ring in $n$
indeterminates over $R$. If we suppose that $T$ is an $m$-SG ring,
it is easy, by \cite[Theorem 2.1]{Bennis and Mahdou3}, to see that
$n\leq m$.

\begin{Theorem}\label{R[X]1}
If $R[X_1,X_2,...,X_n]$ is an $m$-SG ring then $R$ is an
$(m-n)$-SG ring.
\end{Theorem}
\proof By induction on $n$ is suffices to prove the result for
$n=1$. So, suppose that $R[X]$ is an $m$-SG ring. Let $M$ be an
arbitrary $R$-module. For the $R[X]$-module $M[X]:=M\otimes_RR[X]$
there is an exact sequence of $R[X]$-modules $0\longrightarrow
M[X]\longrightarrow P \longrightarrow M[X]\longrightarrow 0$ where
$\pd_{R[X]}(P)\leq m$. Applying $-\otimes_{R[X]}R$ to the short
exact sequence above and seeing that
$M\cong_RM[X]\otimes_{R[X]}R$, we obtain a short exact
 sequence of $R$-modules with the form $0\longrightarrow M \longrightarrow P\otimes_{R[X]}R\longrightarrow M \longrightarrow 0$ (see that $R$ is a projective
 $R[X]$-module). Moreover, $\pd_R(P\otimes_{R[X]}R)\leq
 \pd_{R[X]}(P)<\infty$. On the other hand, by \cite[Theorem 2.1]{Bennis and
 Mahdou3},
 $\Ggldim(R)=\Ggldim(R[X])-1\leq m-1$. Hence, by Proposition
 \ref{prop1}, $R$ is an $(m-1)$-SG ring.\\

 Trivial examples of  $n$-SG-ring  are the
rings with global dimension   $\leq n$. The following example
gives a new family of  $n$-SG rings  with infinite weak global
dimension.
\begin{Example}\label{1}
Consider the non semi-simple quasi-Frobenius rings
$R_1:=K[X]/(X^2)$ and $R_2:=K[X]/(X^3)$ where $K$ is a field, and
let $S$ be a non Noetherian  ring such that $gldim(S)= n$. Then,
\begin{enumerate}
\item $\Ggldim(R_1)=\Ggldim(R_2)=0$ and $R_1$ is $0$-SG ring but
$R_2$ is not.
    \item $R_1\times S$ is a non Noetherian $n$-SG   ring with infinite weak
global dimension.
 \item $\Ggldim(R_2\times S)=n$ but $R_2\times S$ is not an $n$-SG
 ring.
\end{enumerate}

\end{Example}
\proof From \cite[Corollary 3.9]{ouarghi} and \cite[Proposition
2.6]{Bennis and Mahdou2}, $\Ggldim(R_1)=\Ggldim(R_2)=0$ and $R_1$
is $0$-SG ring but $R_2 $ is not. So, $(1)$ is clear.  Moreover
$R_1$ and $R_2$ have infinite weak global dimensions. By
\cite[Theorems 2.1]{Bennis and Mahdou3} and Theorem \ref{product},
it is easy to see that, $\Ggldim(R_2\times S)=n$ and that
$R_2\times S$ is not an $n$-SG
 ring.\\


\begin{thebibliography}{99}



\bibitem{Aus bri}{ M. Auslander and M. Bridger,  Stable module
theory, \em  Memoirs. Amer. Math. Soc., }{\bf 94},  American
Mathematical Society, Providence, R.I., 1969.
\bibitem{Bennis and Mahdou1}{ D. Bennis and N. Mahdou,  Strongly Gorenstein projective,
injective, and flat modules, \em  J. Pure Appl. Algebra}
\textbf{210} (2007), 437--445.
\bibitem{Bennis and Mahdou2}{ D. Bennis and N. Mahdou, Global
Gorenstein Dimensions, \em Proc. Amer. Math. Soc., }\textbf{138}
(2) (2010), 461-465.

\bibitem{Bennis and Mahdou3}{D. Bennis and N. Mahdou, Global Gorenstein dimensions of polynomial rings
and of direct products of rings, \em Houston Journal of
Mathematics, }\textbf{25} (4) (2009), 1019--1028.

\bibitem{ouarghi}{ D. Bennis, N. Mahdou and K.
Ouarghi,  Rings over which all modules are strongly Gorenstein
projective, \em Rocky Mountain Journal of Mathematics, } {\bf 40}
(3) (2010), 749--759.


\bibitem{Christensen}{ L. W. Christensen, \em Gorenstein dimensions, } Lecture Notes
in Math., 1747, Springer, Berlin, 2000.

\bibitem{Christensen and Frankild}{
L. W. Christensen, A. Frankild, and H. Holm, On Gorenstein
projective, injective and flat dimensions - a functorial
description with applications, \em J. Algebra} {\bf 302} (2006),
231–-279.

\bibitem{Enochs}{
E. Enochs and O. Jenda, On Gorenstein injective modules, \em Comm.
Algebra, } {\bf 21} (10) (1993), 3489--3501.

\bibitem{Enochs2}{
E. Enochs and O. Jenda, Gorenstein injective and projective
modules, \em Math. Z., } {\bf 220} (4) (1995),  611--633.
\bibitem{Enocks and janda}{ E. E. Enochs and O. M. G. Jenda, \em Relative Homological Algebra}, de Gruyter Expositions in
Mathematics, Walter de Gruyter and Co., Berlin, 2000.

\bibitem{Eno Jenda Torrecillas}{
E. Enochs, O. Jenda and B. Torrecillas, Gorenstein flat modules,
 \em Nanjing Daxue Xuebao Shuxue Bannian Kan, } {\bf 10} (1)  (1993),  1--9.



\bibitem{Holm}{
H. Holm, Gorenstein homological dimensions, \em J. Pure Appl.
Algebra, } {\bf 189} (2004), 167-–193.

\bibitem{Mahdou costa} {N. Mahdou, On Costa's conjecture, \em Comm. Algebra, }{\bf  29} (7)
(2001), 2775-2785.

\bibitem{MT1}{ N. Mahdou and M. Tamekkante, On (strongly) Gorenstein (semi)hereditary
rings, \em Arab J Sci Eng, }  (Springer) 36 (2011) 431--440.

\bibitem{MahdouTamekkante}{ N. Mahdou and M. Tamekkante, Storongly n-Gorenstein projective, injective and flat modules, \em submitted for publication}. Available from
math.AC/0904.4013 v1 26 Apr 2009.


\end{thebibliography}
\end{document}